\documentstyle[11pt,amscd]{amsart}
\newtheorem{pr}{Proposition}

\newcommand{\proj}{\bold P}
\newcommand{\barr}{\overline}
\newcommand{\rarr}{\rightarrow}
\newcommand{\com}{\Bbb{C}}
\begin{document}
\title{A Geometric Invariant Theory Compactification Of ${ M}_{g,n}$
Via The Fulton-MacPherson Configuration Space}
\author{Rahul Pandharipande$^1$}
\date{10 May 1995}
\maketitle
\pagestyle{plain}
\footnotetext[1]{Research supported by an NSF Graduate Fellowship.}

\setcounter{section}{-1}
\section{Introduction}
\label{int}
In [F-M], a compactification of the
configuration space of $n$ marked points on an algebraic variety
is defined.  For a nonsingular curve $X$ of genus $g \geq 2$,
 the Fulton-MacPherson
configuration space of $X$
(quotiented by the automorphism group of $X$) is isomorphic to
the (reduced)
fiber of $\gamma:\barr{M}_{g,n} \rightarrow \barr{M}_{g}$ over $[X]\in M_g$.
 Since the Fulton-MacPherson configuration
space is defined for singular varieties, it is natural to ask
whether a compactification of $\gamma^{-1}(M_g)$ can be obtained over
$\barr{M}_{g}$.
First, we consider the Fulton-MacPherson configuration
space for families of varieties. This relative construction is then
applied to the universal curve over the Hilbert scheme of
$10$-canonical, genus $g\geq 2$ curves.
Following results of Gieseker, it is shown there
exist linearizations of the natural ${SL}$-action on the relative
configuration space of the universal curve
 that yield G.I.T. quotients compactifying
$\gamma^{-1}(M_g)$. These new compactifications, $M_{g,n}^{c}$, are described.
For $n=1$, $M_{g,1}^c$ and
the Deligne-Mumford compactification $\barr{M}_{g,1}$ coincide. For $n=2$,
$M_{g,2}^c$ and $\barr{M}_{g,2}$ are isomorphic on open
sets with codimension $2$ complements. $M_{g,2}^c$ and $\barr{M}_{g,2}$
differ essentially by the birational modification corresponding
to the two minimal resolutions of an ordinary threefold double point.
For higher $n$, the compactifications
$M_{g,n}^c$ and $\barr{M}_{g,n}$ differ more substantially.

Thanks are due to
J. Harris for mathematical guidance.
The author has benefited from many discussions with him.

\section{Relative Fulton-MacPherson Configuration Spaces}
\label{rfm}

\subsection{Terminology}
\label{rmft}
Let $\com$ be the ground field of complex numbers. A
morphism $\mu: X \rightarrow Y$ is an {\it immersion} if
$\mu$ is an isomorphism  of $X$ onto an open subset of
a closed subvariety of $Y$. A morphism $\gamma$ is {\it quasi-projective}
if it factors as $\gamma =\rho \circ \mu$ where $\mu$ is an open immersion
and $\rho$ is projective.  The only smooth morphisms considered will be
smooth morphisms of relative dimension $k$ between nonsingular
varieties.

\subsection{Definition}
\label{rfmd}
We carry out the construction of Fulton and MacPherson
in  the relative context.
Suppose $\pi:{ F}\rightarrow{ B}$ is  a (separated) morphism
of algebraic varieties. Let $n$ be a positive integer.
${ N}= \{ 1,\ldots, n \}  $. Wherever  possible,
products will be taken in the category of varieties over $B$.
Define:
 $$ F_B^N=\prod_{N}F=
 \underbrace { F \times {_B} F \times {_B} \ldots \times {_B} F}_{n}\;\;. $$
And define:
$$(F_B^N)_0= F_B^N \;\;\setminus \;\; (\bigcup \bigtriangleup_{\{a,b\}})$$
Where $\bigtriangleup_{\{a,b\}}$ denotes the large diagonal
corresponding to the indices $a,b \in N$ and the union is over
all pairs $\{ a,b \}$ of distinct element of $N$.
For each subset $S$  of $N$ define
$F_B^S=\prod_{S}F$. Following the notation of [F-M],
let $Bl_\bigtriangleup (F_B^S)$ denote the blow-up of $F_B^S$ along
the small diagonal. There exists a natural immersion:
\begin{equation}
(F_B^N)_0 \subset F_B^N \times \prod_{|S|\geq 2}{Bl_\bigtriangleup
(F_B^S)}\;\; .\label{mouse}
\end{equation}
The relative Fulton-MacPherson configuration space of
$n$ marked points of $F$ over $B$, $F_B[n]$, is defined
to be the closure of $(F_B^N)_0$ in the above product.  When $B$ is
a point, this definition coincides with that of [F-M].
Consider the composition:
\begin{equation}
F_B[n] \subset  F_B^N \times \prod_{|S|\geq 2}{Bl_\bigtriangleup
(F_B^S)}
\stackrel{\mu}{\rarr}
 F_B^N \times \prod_{|S|\geq 2}{(F_B^S)}
\stackrel{\beta}{\rarr}
 F_B^N
\end{equation}
where $\mu$ is the natural blow-down morphism and $\beta$ is
the projection on the first factor. Since $\mu$ is a projective
morphism and $F_B[n]$ is a closed subvariety,
$$\mu: F_B[n] \rarr \mu(F_B[n])$$
is also projective. Since $\beta: \mu(F_B[n]) \stackrel{\sim}
{\rightarrow} F_B^N$ is an isomorphism, the morphism $\rho=\beta
 \circ \mu$
$$\rho: F_B[n] \rarr F_B^N$$
is projective.
 For our
purposes, we shall only consider the case where
$\pi:{F}\rightarrow{ B}$ is a
quasi-projective morphism. Also, we will
be mainly interested in the case where
$F$ and  $F^N_B$ are irreducible varieties.

\subsection{The Blow-Up Construction}
\label{rfmb}
Consider again the birational projective morphism
$$\rho :F_B[n]\rightarrow F_B^N$$
It is natural to inquire
whether $\rho$ can be expressed as a composition of explicit
blow-ups along canonical subvarieties. In [F-M], such a blow-up
construction is given for the configuration space in case
$B$ is a point. The blow-ups in [F-M] are canonical in the following
sense: if $Y \rightarrow X$ is an immersion, the sequence of
blow-ups resolving $Y[n] \rightarrow Y^N$ is the sequence of
strict transformations of $Y^N$ under the blow-ups resolving
$X[n] \rightarrow X^N$.

The blow-up construction of Fulton and MacPherson is valid in the
relative context.
We now assume that $\pi: F\rightarrow B$ is a quasi-projective morphism.
In this case, there exists a factorization:
\begin{equation*}
\begin{CD}
F @>i>> {\proj}^r \times B \\
@VV{\pi}V       @VVV \\
B @= B
\end{CD}
\end{equation*}
where $i$ is an immersion.
We use the notation $\proj^r \times B = \proj_B^r$ and drop
extra $B$ subscripts when the meaning is clear. For example,
$(\proj_B^r)^N$ instead of $(\proj_B^r)_B^N$ . We have
the following commutative diagram:
\begin{equation}
\begin{CD}
(F_B^N)_0 @>>> F_B^N \times \prod_{|S|\geq 2}{Bl_\bigtriangleup (F_B^S)}\\
@VV{i^N}V             @VV{i^{Bl}}V \\
(\proj_B^r)^N_0 @>>> (\proj_B^r)^N \times \prod_{|S|\geq 2}{Bl_\bigtriangleup
((\proj_B^r)^S)}
\label{snake}
\end{CD}
\end{equation}
where $i^N$, $i^{Bl}$ are immersions. We conclude from
diagram (\ref{snake}) that $F_B[n]$ is immersed in $\proj_B^r[n]$.
Hence:
\begin{equation}
\begin{CD}
F_B[n] @>j>> \proj_B^r[n] \\\
@VV{\rho}V   @VV{\eta}V \\
F_B^N @>{i^N}>> (\proj_B^r)^N
\label{eagle}
\end{CD}
\end{equation}
where $i^N$, $\;j$ are immersions. Since $\rho$ is
a projective morphism, $j(F_B[n])$ is closed
in $\eta^{-1}(i^N(F_B^N))$.
$F_B[n]$ is therefore the strict transformation of $F_B^N$ under $\eta$.
It is clear the following diagram holds:
\begin{equation*}
\begin{CD}
\proj_B^r[n] @= \proj^r[n] \times B \\
@VV{\eta}V    @VV{\gamma \times id}V \\
(\proj_B^r)^N @=( \proj^r)^N \times B
\end{CD}
\end{equation*}
In [F-M], an explicit and canonical blow-up construction of $\gamma$
is given. By extending each exceptional locus over the base
$B$, a blow-up construction of $\eta$ is obtained.  We see from
diagram (\ref{eagle}) that a blow-up construction of $F_B[n]$ exists by
taking the strict transformation of $F_B^N$ at each blow-up of
$(\proj_B^r)^N$.

\subsection{Comparing $F_b[n]$ and $F_B[n]_b$}
\label{rfmc}
For a given $b \in B$ let $F_b$ denote the fiber of $\pi$ over $b$.
{}From equation (\ref{mouse}) and the definitions, it is clear there exists
a natural closed immersion:
$$F_b[n] \stackrel{i_b}{\hookrightarrow}  F_B[n]_b .$$
It is possible for $i_b$ to be a proper inclusion.  Examples of this
behavior will be seen in section (\ref{disc}).
 We have the following:
\begin{pr}
\label{pfiber}
If $B$ is irreducible and $\pi:F \rightarrow B$ is a smooth,
quasi-projective morphism of nonsingular varieties, then
for every $b \in B$, $i_b$ is an isomorphism.
\end{pr}
\begin{pf}
Suppose $X$ is a fixed nonsingular algebraic variety. In [F-M],
the canonical construction of $X[n]$ is given by
 a sequence of explicit blow-ups
of $X^N$ along {\it nonsingular} centers.  In the previous section, it was
shown how the construction of [F-M] could be lifted to the relative
context. Let $m$ be the number of blow-ups needed in the Fulton-
MacPherson construction resolving
$\rho: F_B[n] \rightarrow F_B^N$.  Let $F_{B,j}^N$ for
$0 \leq j \leq m$ denote the $j^{th}$ stage. $F_{B,0}^N= F_B^N$
and $F_{B,m}^N= F_B[n]$. Since the blow-up construction in [F-M] is
canonical, for any variety $X$ similar definitions can be made. We
show inductively, for each $b \in B$, the natural inclusion:
\begin{equation}
\label{duck}
F_{b,j}^N \hookrightarrow (F_{B,j}^N)_b
\end{equation}
is an isomorphism.  For $j=0$ the assertion is clear.  The induction
step rests on a simple {\bf Claim}:

Suppose $S$ is an irreducible nonsingular variety,
 $ R \rightarrow S$ is a smooth morphism,
and $T \hookrightarrow R$ is a closed immersion smooth over $S$ . Then, for
any $s \in S$, the blow-up of $R_s$ along $T_s$ is naturally isomorphic to
the fiber over $s$ of the blow-up of $R$ along $T$. Since all spaces are
nonsingular, the assertion follows from examining normal directions of $T$
in $R$; the various smoothnesses imply all normal directions are
represented in the fiber.

Assume equation (\ref{duck}) is an isomorphism for all $b \in B$.
Since $\pi$ is smooth of relative dimension $k$,
 $F_b$  and $F_{b,j}^N$ are nonsingular of pure dimensions $k$, $nk$. Hence,
$(F_{B,j}^N)_b$ is nonsingular of pure dimension $nk$.
Also, every irreducible
component of $F_B^N$ (and hence $F_{B,j}^N$) is of relative
dimension $nk$. The last two facts
   imply the morphism: $$\pi_j^N:F_{B,j}^N
\rightarrow B$$ is smooth. Examination of the $(j+1)^{th}$ center is
straightforward.  Because of the assumed isomorphism (\ref{duck}) and the
knowledge that the [F-M] construction of the configuration space of a
nonsingular variety over a point only involves nonsingular centers, we see
that the $(j+1)^{th}$ center is smooth over B.  The above claim now proves
the induction step.
\end{pf}

\subsection{Universal Families}
\label{rfmu}
Let $X$ be a nonsingular algebraic variety. Let
$\barr{x}=(x_1, \ldots, x_n)$ be  $n$ ordered points of $X$.
A subset $S\subset N$ is said to be coincident
 at $z\in X$ if $|S| \geq 2$ and
for all $i \in S$, $x_i=z$. Following [F-M], for every $S$ coincident
at $z$, we define a {\it screen} of $S$ at $z$ to be an
equivalence class of the data $(t_i)_{i\in S}$ where:
\begin{enumerate}
\item  $t_i \in T_z$, the tangent space of $X$ at $z$.
\item  $\exists i,j \in S$ such that $t_i \neq t_j$.
\end{enumerate}
Two data sets $(t_i)_{i\in S}$ and $(t'_i)_{i\in S}$
are equivalent if there exists $\lambda \in C^*$ and
$v \in T_z$ so that for all $i \in S$,
$\lambda \cdot t_i + v= t'_i$.
A screen shows the tangential separation of infinitely
near points.
An $n$-tuple $\barr{x}=(x_1, \ldots, x_n)$ together
with  a screen $Q_S$ for each  coincident subset  $S\subset N$
constitute an {\it n-pointed stable class} in $X$ if the
screens are compatible in the following sense.
Suppose $S_1\subset S_2$ are two subsets coincident at
$z$ where $Q_{S_2}$ is represented by the data  $(t_j)_{j \in S_2}$.
If there exist
$k,\hat{k} \in S_1$ so that $t_k\neq t_{\hat k}$,
then $(t_j)_{j\in S_1}$ defines a screen for $S_1$.
The compatibility condition requires that when
this restriction of $Q_{S_2}$ is defined, it equals  $Q_{S_1}$.
For a nonsingular space $X$, $X[n]$ is the parameter
space of $n$-pointed stable classes in $X$.
Given an $n$-pointed stable class in $X$, an
{\it $n$-pointed stable degeneration } of $X$ can be constructed
(up to isomorphism)
 as follows. Let $z\in X$ occur with multiplicity in $\barr{x}$.
Blow-up $X$ at $z$ and attach a
 $\proj(T_z \oplus {\bold 1})$ in the natural way along
the exceptional divisor at $z$.  Note that
$\proj(T_z \oplus {\bold 1})$ minus the hyperplane at infinity,
$\proj(T_z)$, is naturally isomorphic to the affine space
$T_z$. Let $S_z \subset N$ be the maximal subset coincident
at $z$. The screen $Q_{S_z}$ associates (up to equivalence)
points of $T_z$ to the indices that lie in $S_z$. Condition
(2) of the screen data implies some separation of the
marked points has occurred. The further screens specify
in a natural way (up to equivalence of screens)
the further blow-ups and markings required to separate the marked points.
The final space obtained along with $n$ distinct marked points is
the $n$-pointed stable degeneration associated to the given
$n$-pointed stable class.
See [F-M] for further details.

It is shown in [F-M] that if $X$ is an nonsingular
variety, there exists a universal family of $n$-pointed
stable degenerations of $X$ over $X[n]$. Let $X[n]^+$ denote
this universal family. $X[n]^+$ is equipped with the following
maps:
\begin{equation*}
\begin{CD}
X[n]^+ @>{\mu}>> X[n] \times X \\
@VV{\mu_p}V     @VVV \\
X[n] @= X[n]
\end{CD}
\end{equation*}
There are $n$ sections of ${\mu_p}$, $\{ \sigma_i \}_{i \in N}$ :
$$X[n] \stackrel{\sigma_i}{\rightarrow} X[n]^+ . $$
For any $d\in X[n]$, the fiber $\mu^{-1}_p(d)$ along with
the $n$-tuple $(\sigma_1(d), \ldots,\sigma_n(d))$ is
the $n$-pointed stable degeneration of $X$ associated to
the $n$-pointed stable class corresponding to $d$.

 We note here that if $C$ is
a nonsingular automorphism-free curve,   $n$-pointed stable classes
in $C$ correspond bijectively to isomorphism classes of
 $n$-pointed Deligne-Mumford
stable curves over $C$ . Moreover, the universal family over $C[n]$ defines
a map to the reduced fiber $$\phi:C[n] \rightarrow \gamma^{-1}([C])$$ where
$\gamma :\barr{M}_{g,n} \rightarrow \barr{M}_g$.
Since $\phi$  is proper bijective
and both spaces are normal, $\phi$ is an isomorphism. If $C$
has a finite automorphism group, $A$, we see $A$ acts on $C[n]$
and $\phi$ is $A$-invariant. Therefore $\phi$ descends to the quotient
$$\phi (C[n]/A) \rightarrow \gamma^{-1}([C]).$$
It is not hard to see that this map is proper bijective and hence
an isomorphism by  normality.

The map $\mu$ is a birational morphism and can be expressed as an explicit
sequence of blow-ups of $X[n] \times X$ along canonical, nonsingular loci.
Canonical here has the same meaning as in section (\ref{rfmb}) :
if $Y \rightarrow X$ is an immersion of nonsingular varieties, the
blow-up sequence resolving $Y[n]^+ \rightarrow Y[n]\times Y$ is
the strict transform of of $Y[n] \times Y$ under the blow-up
sequence resolving $X[n]^+ \rightarrow X[n] \times X$. Moreover,
the sections of $Y[n]^+ \rightarrow Y[n]$ are restrictions
of the sections of $X[n]^+ \rightarrow X[n]$. This
canonical blow-up construction is given in [F-M].

\subsection{Relative Universal Families}
\label{rfmru}
Suppose $\pi:F \rightarrow B$ is a smooth,
quasi-projective morphism of  nonsingular varieties. In this case,
the construction of the
universal family that appears in [F-M] can be lifted to the relative
context. Using the notation of section (\ref{rfmb}), we have an immersion:
$$F_B[n] \times_B F \rightarrow \proj_B^r[n] \times_B \proj_B^r \;.$$
Consider the diagram:
\begin{equation*}
\begin{CD}
\proj_B^r[n]^+ @= \proj^r[n]^+ \times B \\
@VV{\omega}V    @VV{\mu \times id}V \\
\proj_B^r[n] \times_B \proj_B^r  @= \proj^r[n]\times \proj^r
 \times B
\end{CD}
\end{equation*}
For $\omega=\mu \times id$, the Fulton-MacPherson construction of the
universal family can be carried out uniformly over the base by
extending the centers of the blow-ups resolving $\mu$
 trivially over $B$.  Define
$F_B[n]^+$ to be the
proper transform of $ F_B[n] \times_B F $ under $\omega$. We have:
$$\upsilon: F_B[n]^+ \rightarrow F_B[n] \times_B F$$
To show the space defined above, $F_B[n]^+$, has the
desired geometrical properties, we argue as in the proof
of Proposition 1. Let $ (F_B[n] \times_B F)_j $ denote the
$j^{th}$ stage of the canonical sequence of blow-ups resolving $\upsilon$.
Inductively, it is shown that for each $b \in B$ there is an isomorphism:
\begin{equation}
\label{cow}
 (F_b[n] \times F_b)_j \rightarrow
 (F_B[n] \times_B F)_{j,b}.
\end{equation}
The $j=0$ case is established by Proposition 1. The induction step follows
from the the claim made in the proof of Proposition 1 and the fact
that for a nonsingular variety $X$, the canonical Fulton-MacPherson resolution
of $X[n]^+ \rightarrow X[n]\times X$ involves only nonsingular centers.

We conclude that  fiber  $F_B[n]^+_b$ over $b \in B$ is
naturally  isomorphic
to $F_b[n]^+$.
It is clear that $n$ sections ${\sigma_i}$ exist for
 $$\omega_p: \proj_B^r[n]^+ \rightarrow \proj_B^r[n].$$
For each $b \in B$, these sections
${\sigma_i}$ are compatible with the $n$ natural sections of
$F_b[n]^+ \rightarrow F_b[n].$
 Therefore, via restriction, the ${\sigma_i}$ yield $n$ sections of
 $$\upsilon_p: F_B[n]^+ \rightarrow F_B[n].$$
The fiber of $F_B[n]^+_\xi$ over
$\xi \in F_B[n]$ is a $n$-pointed stable degeneration of
$F_{\pi(\xi)}$. We have:
\begin{pr}
Suppose $B$ is irreducible and $\pi: F \rightarrow B$ is a smooth,
quasi-projective  morphism of nonsingular varieties, then
 $F_B[n]^+$ along with $\upsilon$ and
 $\{\sigma_i\}_{i \in N}$ is a universal family of n-pointed stable
degenerations of $F_B$ over $F_B[n]$.
\end{pr}

\subsection{Final Note}
\label{rfmfn}
Suppose $\pi:G \rightarrow B$ is a projective morphism,
$G$ is nonsingular, irreducible, $B$ is nonsingular, $\pi$ is flat, and
for every $b \in B$ the fiber $G_b$ is reduced.
In this case, let $F \subset G$ be the open set where
$\pi$ is smooth. Using flatness and a tangent space
calculation, we see:
$$ F=\{\xi \in G|\xi \hbox{ is a nonsingular point of } G_{\pi(\xi)} \} $$
and  $\pi: F \rightarrow B$ is a smooth, surjective
morphism of nonsingular varieties.
We know the space $F_B[n]$ is equipped with a
universal family $F_B[n]^+$ obtained from $F_B[n]\times_B F$
by a sequence of canonical blow-ups.  The problem with this
universal family is that its fibers over $F_B[n]$ are $n$-pointed
stable degenerations of $F_B$ not $G_B$.  This problem can easily
be fixed. Note there is an open inclusion:
$$F_B[n]\times_B F \subset F_B[n]\times_B G .$$
It is the case that the centers of the blow-ups resolving
$$\upsilon:F_B[n]^+ \rightarrow F_B[n] \times_B F$$
are closed in $F_B[n]\times_B G$.
Using the isomorphism  (\ref{cow}) and the explicit
description of the centers of blow-ups in [F-M], this
closure is not hard to check.
Hence, if the sequence of
blow-ups is carried out over $F_B[n]\times_B G$ the desired
family of $n$-pointed stable degenerations of $G_B$ is obtained
over $F_B[n]$. An $n$-pointed stable degeneration of a fiber
$G_b$ is as before with the additional condition
that  the marked points must lie
 over the smooth locus of $G_b$.

\section{The Geometric Invariant Theory Set-Up}
\label{git}

\subsection{Notation}
\label{gitn}
 Let ${\barr{ M}_g}$
denote the Deligne-Mumford compactification of the moduli space of
nonsingular, genus $g$, projective  curves, $M_g$. Let ${\barr{M}_{g,n}}$
denote the Deligne-Mumford compactification of the moduli space of genus
$g$ curves with $n$ marked points. There exists a natural projective
morphism
$$\gamma : {\barr {M}_{g,n}} \rightarrow {\barr {M}_g}.$$
All these spaces are normal.

\subsection{Gieseker's construction of $\overline{M}_g$}
\label{gitg}
Fix an integer $g \geq 2$. Define: $$d=10 \cdot (2g-2)$$
$$R=d-g.$$ Define the polynomial:
$$f(m)= d \cdot m -g +1.$$
Note $f(m)$ is the Hilbert polynomial of a complete, genus $g$,
 $10$-canonical curve in $\proj ^R$. Let $H_{f,R}$ denote
the Hilbert scheme of the polynomial $f$ in $\proj ^R$.
If $X$ is a closed subscheme of $\proj ^R$ with Hilbert polynomial
$f$, we denote the point of $H_{f,R}$ corresponding to $X$ by $[X]$.
It is well known that there exists an integer ${\widehat{m}}$ such
that, for any $m \geq \widehat{m}$ and any closed subscheme
 $X \subset \proj ^R$ corresponding to a point  $ [X] \in H_{f,R}$,
\begin{equation}
\label{hippo}
h^1(I_X(m),\proj ^R) = 0
\end{equation}
\begin{equation}
\label{bear}
h^0({\cal O}_X(m), X)= f(m) \; .
\end{equation}
Therefore, for any $m \geq \widehat{m}$, there is a natural map:
$$i_m: H_{f,R} \rightarrow
 {\proj}(\bigwedge ^{f(m)} H^0( {\cal O}_{{\proj}^R}(m), \proj ^R)^*).$$
Where $i_m$ is defined  for each $[X] \in H_{f,R}$ as follows:
by (\ref{hippo}), there is a natural surjection
$$ H^0({\cal O}_{\proj ^R}(m), \proj ^R) \rightarrow
H^0({\cal O}_X(m),X) $$
which yields, by (\ref{bear}), a surjection
\begin{equation}
\label{whale}
 \bigwedge ^{f(m)} H^0({\cal O}_{\proj ^R}(m), \proj ^R) \rightarrow
\bigwedge ^{f(m)} H^0({\cal O}_X(m),X) \cong {\bold C}.
\end{equation}
The last surjection (\ref{whale}) is an element of
${\proj}(\bigwedge ^{f(m)} H^0( {\cal O}_{{\proj}^R}(m), {\proj}^R)^*)$.
The map $i_m$ is now defined on sets. That $i_m$ is an algebraic morphism
of schemes can be seen by constructing (\ref{whale}) uniformly
over $ H_{f,R}$ and  using the universal property of
${\proj}(\bigwedge ^{f(m)} H^0({\cal O}_{{\proj}^R}(m), {\proj}^R)^*)$.
In fact, it can be shown there
 exists an integer $\overline{m}$ such that for every
$m \geq \overline{m}$, $i_m$ is a closed immersion.

{}From the universal property of the Hilbert scheme, we obtain
a natural $SL_{R+1}$-action on $H_{f,R}$.
 For each $m \geq \overline{m}$,
the closed immersion $i_m$ defines a linearization of the natural
$SL_{R+1}$-action on $H_{f,R}$.
Define the following locus $\barr{K}_g\subset H_{f,R}$:
$[X] \in \barr{K}_g$ if and only if $X$  is a nondegenerate, 10-canonical,
genus g, Deligne-Mumford stable curve in $\proj^R$.
$\barr{K}_g$ is a quasi-projective, $SL_{R+1}$-invariant subset. In [G],
Gieseker shows a linearization $i_m$ can be chosen satisfying:
\begin{enumerate}
\item[(i)] The stable locus of the corresponding G.I.T. quotient contains
$\barr{K}_g$.
\item[(ii)] $\barr{K}_g$ is closed in the semistable locus.
\end{enumerate}
{}From (i), we see $\barr{K}_g/SL_{R+1}$ is a geometric quotient.
By (ii), $\barr{K}_g/SL_{R+1}$ is a projective variety. Since
$\barr{K}_g$ is a nonsingular variety ([G]), it follows that
  $\barr{K}_g/SL_{R+1}$ is normal.
{}From the definition
of $\barr{K}_g$, the universal family over $H_{f,R}$
 restricted to $\barr{K}_g$ is a
family of Deligne-Mumford stable curves.  Therefore there exists a
natural map $\mu: \barr{K}_g \rightarrow \overline{M}_g$. Since $\mu$ is
$SL$-invariant,  $\mu$ descends to a projective morphism
from the quotient $\barr{K}_g/SL_{R+1}$ to $\overline{M}_g$. Since $\mu$
is one to one and $\overline{M}_g$ is normal, $\mu$ is an
isomorphism. Note that since $\barr{M}_g$ is irreducible,
$\barr{K}_g$ is also irreducible.

\subsection{The Relative $n$-pointed Fulton-MacPherson Configuration
Space of the Universal Curve}
\label{gitrf}
Let $\pi: U_{H} \rightarrow H_{f,R}$  be the universal
family over the Hilbert scheme defined in section (\ref{gitg}) where
$\pi$ is a flat, projective morphism.
 Let $\barr{K}_g \subset H_{f,R}$
be defined as above. Let $U_{\barr{K}_g}$ be the
 restriction of $U_H$ to $\barr{K}_g$. Following the notation of
section (\ref{rfmd}), we define $U_{\barr {K}_g}[n]$ to be the relative Fulton-
MacPherson space of $n$-marked points on $U_{\barr {K}_g}$ over
 ${\barr {K}_g}$.
{}From section (\ref{rfmb}), we see the immersion $\zeta$:
\begin{equation*}
\begin{CD}
U_{\barr {K}_g} @>{\zeta}>> \proj ^R \times H_{f,R} \\
@VV{\pi}V  @VV{\rho}V \\
{\barr {K}_g} @>>> H_{f,R}
\end{CD}
\end{equation*}
yields another immersion $\zeta[n]$:

\begin{equation*}
\begin{CD}
U_{\barr{K}_g}[n] @>{\zeta[n]}>> \proj ^R[n] \times H_{f,R} \\
@VV{\pi[n]}V  @VV{\rho[n]}V \\
{\barr{K}_g} @>>> H_{f,R}
\end{CD}
\end{equation*}
There exists a natural $SL_{R+1}$-action  on  $ \proj ^R[n]$
and therefore on $ \proj ^R[n] \times H_{f,R}$. Since $U_{\barr {K}_g}$ is
 invariant under the natural $SL_{R+1}$-action, we see $U_{\barr{K}_g}[n]$
is also $SL_{R+1}$-invariant.
Since $\pi$ is projective,
 $U_{\barr{K}_g}[n] \subset \rho[n]^{-1}(\barr {K}_g)$ is
a closed subset. It follows from (i) and (ii) of section (\ref{gitg})
and Propositions (7.1.1) and (7.1.2) of [P] that there exist linearizations
of the natural $SL_{R+1}$-action on $ \proj ^R[n] \times H_{f,R}$ satisfying:
\begin{enumerate}
\item[(i)] $U_{\barr{K}_g}[n]$ is contained in the stable locus of the
corresponding G.I.T. quotient.
\item[(ii)] $(\rho[n]^{-1}(\barr{K}_g))^{SS}$
 is closed in the semistable locus.
\end{enumerate}
{}From (i), (ii), and the fact that $U_{\barr{K}_g}[n]$ is closed in
 $\rho[n]^{-1}(\barr {K}_g)$,
we see that $U_{\barr {K}_g}[n]/SL_{R+1}$ is a geometric quotient
and a projective variety. Define:
$$M_{g,n}^{c} = U_{\barr {K}_g}[n]/SL_{R+1} .$$
Note there is a natural projective morphism
$$\rho : M_{g,n}^{c} \rightarrow {\barr {M}_{g}}$$
descending from the $SL_{R+1}$-invariant maps:
$$U_{\barr {K}_g}[n] \rightarrow {\barr {K}_g}
\rightarrow {\barr {M}_g}.$$
It follows easily that $M_{g,n}^{c}$ is a compactification
of $\gamma^{-1}(M_g)$.  To see this first make the definition:
$${K_g} =\{[X] \in H_{f,R}|X \hbox{ is a nondegenerate, 10-canonical,
nonsingular, genus g curve}  \}. $$
 $U_{K_g}[n]$ is a dense open $SL_{R+1}$-invariant subset of
 $U_{\barr{K}_g}[n]$.
Since the morphism $\pi :U_{K_g} \rightarrow K_g$ is smooth, we see from
section (\ref{rfmru}) that there exists a universal family of
 Deligne-Mumford stable $n$-pointed genus $g$ curves over $ U_{K_g}[n]$.
This universal family yields a canonical morphism
$$\mu: U_{K_g}[n] \rightarrow \gamma^{-1}(M_g).$$
It is easily checked that $\mu$ is $SL_{R+1}$-invariant. Therefore,
$\mu$ descends to the open set, $\rho^{-1}(M_g)$, of $M_{g,n}^{c}$. One
sees $$\mu_d:\rho^{-1}(M_g)  \rightarrow \gamma^{-1}(M_g)$$ is
a bijection by Proposition (\ref{pfiber})
 and the fact  that, for a smooth curve C,
$$ (C[n]/\hbox{automorphisms}) \cong  \gamma^{-1}([C])\subset
\gamma^{-1}(M_g). $$
(See section (\ref{rfmu})).
Since $\rho: \rho^{-1}(M_g) \rightarrow M_g$ is projective,
$\gamma: \gamma^{-1}(M_{g}) \rightarrow M_g$ is separated,
and $\rho= \gamma \circ \mu_d$, we conclude  $\mu_d$
is projective. A bijective projective morphism onto a normal
variety is an isomorphism. Since $\gamma^{-1}(M_g)$ is normal, $\mu_d$ is
an isomorphism.

\section{A Description Of $M_{g,n}^c$}
\label{disc}
\subsection{}
\label{dis1}
Let $\pi: U_{\barr{K}_g} \rightarrow \barr{K}_g$ be as above.
Following  section (\ref{rfmfn}), we define $F\subset U_{\barr{K}_g}$
to be the locus where $\pi$ is smooth.  $F_{\barr{K}_g}[n]
\subset U_{\barr{K}_g}[n]$ is an open $SL$-invariant subset.
The points of $F_{\barr{K}_g}[n]$ parameterized $n$-pointed stable
classes on the nonsingular locus of the fibers of $\pi$. There
exists a universal family over $F_{\barr{K}_g}[n]$ which defines an
$SL$-invariant morphism :
\begin{equation*}
\mu :F_{\barr{K}_g}[n] \rightarrow \barr{ M}_{g,n}^{s}.
\end{equation*}
Where $\barr{ M}_{g,n}^{s}$
parameterizes $n$-pointed, genus g,
Deligne-Mumford stable curves with  marked points lying  over
nonsingular points of the contracted stable model. Let
$$F_{\barr{K}_g}[n]/SL_{R+1} = (M_{g,n}^c)^{s}.$$
$SL$-invariance implies $\mu$ descends to:
\begin{equation}
\label{quail}
\mu_d:  (M_{g,n}^c)^{s} \rightarrow  \barr{ M}_{g,n}^{s}.
\end{equation}
{}From the arguments of section (\ref{rfmu}), we see
$\mu_d$ is bijective. From the valuative criterion,
it follows  $\mu_d$ is proper.
As before, by normality, it follows that $\mu_d$ is
an isomorphism.

\subsection{Points Of $M_{g,n}^c$ Over A Singular Point}
\label{dissing}
{}From section (\ref{dis1}), it is clear only the behavior
of $U_{\barr {K}_g}[n]$ over a singular point of
$U_{\barr{K}_g}$ remains to be investigated. Since this is
a local question about the the smooth deformation
of a node, it suffices to investigate the family:
\begin{equation*}
\begin{CD}
G @>>> Spec(C[x,y]) \times Spec(C[t]) \\
@VV{\pi}V @VVV \\
Spec(C[t]) @= Spec(C[t])
\end{CD}
\end{equation*}
Where $G$ is defined by the equation $xy-t$.
In the Fulton-MacPherson configuration space
$Spec(C[x,y])[n]$, there is a closed subset $\cal {T}_n$
corresponding to the points lying over $(0,0)$.
In the notation of section (\ref{rfmb}),
$$ \cal{T}_n = \rho^{-1}(\underbrace{(0,0),(0,0), \ldots,(0,0)}_{n}).$$
Recall the notation of section (\ref{rfmd}).
Let $B=Spec(C[t])$, $B^*=Spec(C[t])-(0)$, and $G^*=\pi^{-1}(B^*)$.
We want to investigate the subset  $\cal{W}_n \subset (\cal{T}_n,0)$ that
lies in the closure of $G_{B^*}^{*}[n]$ in
$Spec(C[x,y])[n]\times B$.

Suppose $\kappa$ is a family in $(G_{B^*}^{*N})_0$ where
all the marked points specialize to the node $\zeta$ of $G_0$.
After a base change, $t \rightarrow t^r$,
 $\kappa$ can be defined by $n$ sections,
 $(\kappa_1, \ldots, \kappa_n)$,
of $\pi$ in a neighborhood of $0 \in B$.
 The equation of $G$ after base change is now
 $G_r=xy-t^r$.
Let us take $r=2$. The blow-up
of $G_2$ at $\zeta$ is nonsingular and  is
defined in an open set
by the equation $ab-1$ in $Spec(C[a,b]) \times Spec(C[t])$.
The
blow-down morphism is defined by the equations:
$$x=at$$
$$y=bt.$$
Now assume that the strict transforms of the
sections, $(\kappa_1, \ldots, \kappa_n)$, meet
the exceptional curve $(ab=1,t=0)$ in distinct points
$((a_1,a_1^{-1}), \ldots, (a_n,a_n^{-1}))$,
$\forall i$ $a_i\neq 0$. Then it is clear that the
$n$-pointed stable class in $\cal{T}_n$ that is the limit of $\kappa$
is the class in the tangent space of $C[x,y]$ at $(0,0)$
defined by the pairs of vectors:
$$((a_1,a_1^{-1}), \ldots, (a_n,a_n^{-1}))$$
in the basis $(\frac{\partial}{\partial x}, \frac{\partial}{\partial y})$.

We now define a map:
\begin{equation*}
\theta_n: (C^{*N})_0 \rightarrow \cal{T}_n
\end{equation*}
Where $\theta_n((a_1,\ldots,a_n))$ is the $n$-pointed
stable class defined by the tangent vectors
$((a_1,a_1^{-1}), \ldots, (a_n,a_n^{-1}))$.
The preceding paragraph shows that $Image(\theta_n)\subset \cal{W}_n$.
In fact, it is not hard to see that $Image(\theta_n)$
is dense in a component of $\cal{W}_n$.
For $n=2$, $\cal{W}_2= \cal{T}_2$ where $\cal{T}_2$ is just
the $\proj^1$ of normal directions.

Suppose $n \geq 3$. Let $\hat{a}=(a_1,\ldots,a_n)$ and
$\hat{b}=(b_1,\ldots,b_n)$ be distinct points of $(C^{*N})_0$.
Then, $\theta_n(\hat{a})=\theta_n(\hat{b})$ if and only if
there exists a tangent vector $(v_1,v_2)$
and an element  $\lambda \in C^*$  such
that :
$$ \forall i, \;\;\; \lambda \cdot a_i+v_1=b_i \hbox{ and }
\lambda \cdot a_i^{-1}+v_2=b_i^{-1} .$$
These equations imply
\begin{equation}
\label{elk}
\forall i,j,\; \;\;\; \lambda \cdot (a_i-a_j)=(b_i-b_j)
\end{equation}
\begin{equation}
\label{rhino}
\forall i,j,\; \;\;\; \lambda \cdot(a_i^{-1}-a_j^{-1})=(b_i^{-1}-b_j^{-1})
\end{equation}
Dividing (\ref{elk}) by (\ref{rhino}) yields $a_i\cdot a_j=b_i\cdot b_j$.
For $n\geq 3$, we easily obtain $\hat{a}=\pm\hat{b}$.
 Therefore, a component of
$\cal{W}_n$ can  be viewed as  a compactification of $(C^{*N})_0/(\pm)$.
We note that the dimension of $\cal{W}_n$ is $n$ for
$n\geq 3$.

\section{Comparison with $\barr {M}_{g,n}$ for $n=1,2$}
\subsection{$n=1$}
{}From the definitions, $M_{g,1}^c$ equals $U_{\barr {K}_g}/SL_{R+1}$.
$\pi^*(U_{\barr{K}_g})$ is a family of 1-pointed Deligne-Mumford
genus $g$ curves over $ U_{\barr{K}_g}$ via the natural diagonal section.
This tautological family yields an $SL$-invariant morphism:
$$\mu:U_{\barr{K}_g} \rightarrow \barr{M}_{g,1} $$
that descends to
$$\mu_d: M_{g,1}^c \rightarrow \barr{M}_{g,1}.$$
Since $\mu_d$ is proper bijective and $\barr{M}_{g,1}$ is normal,
$\mu_d$ is an isomorphism.

\subsection{$n=2$}
Consider the family:
\begin{equation}
U_{\barr {K}_g}^2=U_{\barr {K}_g} \times_{\barr {K}_g} U_{\barr {K}_g}\;.
\label{mongoose}
\end{equation}
The singular locus of $U_{\barr {K}_g}^2$, $S$, is nonsingular
of pure codimension $3$ and $SL_{R+1}$-invariant. The singular points are pairs
$(\zeta,\zeta)$ where $\zeta \in U_{\barr {K}_g}$ is a node
of a fiber. Moreover, the singularities of $ U_{\barr {K}_g}^2$
are \'etale-locally ordinary threefold double point
singularities. That is, the
singularities are of the form
\begin{equation}
\label{possum}
W\times Spec(C[a,b,c,d]/(ab-cd)) \subset W\times Spec(C[a,b,c,d])
\end{equation}
Where $W$ is nonsingular.  These assertions about the
singular locus follow from the deformation theory of a Deligne-Mumford
stable curve and [G].

There are three standard resolutions of
the ordinary double point  singularity $Spec(C[a,b,c,d]/(ab-cd))$:
\begin{enumerate}
\item The blow-up along  $(a,b,c,d)$.
\item For any $\lambda \in C$, the blow-up along
 $(a-\lambda \cdot c, \lambda \cdot b-d)$.
\item For any $\lambda \in C$, the blow-up along
 $(a-\lambda \cdot d, \lambda \cdot b-c)$.
\end{enumerate}
Methods (2) and (3) yield the distinct small resolutions.
The local description (\ref{possum}) implies that the blow-up
of $U_{\barr {K}_g}^2$ along $S$ is nonsingular with
an exceptional divisor $E$ that is a $\proj^1\times \proj^1$-
bundle over $S$. Using the techniques of
section (\ref{gitrf}), it can be shown that the
natural $SL_{R+1}$-action on the blow-up
$Bl_{(S)}(U_{\barr {K}_g}^2)$ can be linearized so
that all the points in question are stable and the
quotient is projective.
The diagonal embedding
$$D:U_{\barr {K}_g} \hookrightarrow U_{\barr {K}_g}^2$$
is divisorial except along $S$ where it of the form
of (2) and (3) in the local description (\ref{possum}).
By definition,
 $$M_{g,2}^c= Bl_{(D)}(U_{\barr {K}_g}^2)\;\;/SL_{R+1}.$$
There is a natural blow-down map:
 $$\rho: Bl_{(S)}(U_{\barr {K}_g}^2)\rightarrow
Bl_{(D)}(U_{\barr {K}_g}^2).$$
Another $SL_{R+1}$-invariant small resolution of
$U_{\barr {K}_g}^2$ can be obtain by blowing-down uniformly along
the opposite ruling of $E$ blown-down by $\rho$. Let $Y$ denote
this other small resolution and let
$$\barr{\rho}:  Bl_{(S)}(U_{\barr {K}_g}^2) \rightarrow Y$$
be the blow-down.
Linearizations can be chosen so that
$$Y/SL_{R+1} \cong \barr{M}_{g,2}.$$
There are birational morphisms
\begin{equation}
\label{birad}
M_{g,2}^c \ \leftarrow \ \ Bl_{(S)}(U_{\barr {K}_g}^2)\;\;/SL_{R+1}
\ \ \rarr\
\barr{M}_{g,2}.
\end{equation}
Consider the open loci of $M^{c}_{g,2}$ and $\barr{M}_{g,2}$
where the underlying curve has no (nontrivial) automorphism.
On the automorphism free loci the birational modification
(\ref{birad}) is easy to describe. Let $F_1\subset M_{g,2}^c$
be the locus of
of $2$-pointed stable classes that lie over a
node in a  Deligne-Mumford stable curve of genus g. Similarly, let
$F_2\subset \barr{M}_{g,2}$ be the locus of $2$-pointed, genus g,
Deligne-Mumford stable curves such that the
marked points are coincident at a node in
the stable contraction. On the
automorphism free loci,
$M_{g,2}^c$ and $\barr{M}_{g,2}$ are the distinct
small resolution of the fiber product of the universal
curve with itself. Hence, on the automorphism free loci,
the blow-up of $M_{g,2}^2$ along $F_1$
is isomorphic to the blow-up of $\barr{M}_{g,2}$ along
$F_2$.
The modification (\ref{birad}) obtained by this
isomorphism.

\noindent
\address{Department of Mathematics \\ University of Chicago \\
5734 S. University Ave\\ Chicago, IL 60637 \\ rahul@@ math.uchicago.edu}

\end{document}